\overfullrule=0pt
\input xypic
\input amssym
\def\grt{{\frak{grt}}}
\def\ds{{\frak{ds}}}
\def\krv{{\frak{krv}}}
\centerline{\bf The double shuffle Lie algebra injects into the 
Kashiwara-Vergne Lie algebra}
\vskip .4cm
\centerline{Leila Schneps}
\vskip .6cm
\noindent {\bf Abstract.} In this article we prove that
there exists an injective Lie morphism from the double shuffle 
Lie algebra $\ds$ into the Kashiwara-Vergne Lie algebra $\krv$, 
forming a commutative triangle with the known Lie injections of the 
Grothendieck-Teichm\"uller Lie algebra
$\grt\hookrightarrow\ds$ and $\grt\hookrightarrow \krv$\footnote{${}^1$}{A
proof of this injective morphism $\ds\hookrightarrow\krv$ was proposed in
[12], but the proof of the necessary condition
that elements of $\ds$ induce special derivations of ${\rm Lie}[x,y]$
relied on an earlier result by J.~\'Ecalle, published in [5] with 
a sketch of proof.  However no one was able to complete the proof of
\'Ecalle's result, so that the previous proof of the injection 
$\ds\hookrightarrow\krv$ was incomplete. The present text gives a
new, simple proof that elements of $\ds$ induce special 
derivations, not relying on \'Ecalle's earlier result. It turns out 
(personal communication) that \'Ecalle's result has also been proved in 
the uncirculated preprint [6] by Enriquez and Furusho, who use this 
to fill the gap in [12] and therefore obtain the injection 
$\ds\rightarrow \krv$ by the previous method.}.  
\vskip .5cm
\noindent {\bf \S 1. Background and main statement}
\vskip .3cm
\noindent Let ${\rm Lie}[x,y]$ denote the degree-completed free Lie algebra
on two generators $x$ and $y$.
\vskip .2cm
\noindent {\bf Definition of $\ds$.} The {\it double shuffle Lie algebra}
is the space of elements $f\in {\rm Lie}[x,y]$ with no terms of degree $\le 3$, 
satisfying the following property. Let $y_i=x^{i-1}$, so that
every monomial $w$ in $x,y$ ending in $y$ can be written $w=y_{i_1}\cdots 
y_{i_r}$.  Then for any two such words $u,v$, we have
$$\sum_{w\in st(u,v)} (f|w)=0,$$
where the sum is over the stuffle product of the two words $u$ and $v$. The
space $\ds$ is a Lie algebra under the Poisson bracket\footnote{${}^{2}$}{A difficult 
result due independently to \'Ecalle and Racinet, see [8] for Racinet's 
proof and [10] or the mould theory reference text [11] 
for a complete and detailed exposition of 
\'Ecalle's proof.}, defined by
$$\{f,f'\}=D_f(f')-D_{f'}(f)+[f,f']$$
where $D_f$ denotes the derivation of ${\rm Lie}[x,y]$ defined by $D_f(x)=0$,
$D_f(y)=[y,f]$.
\vskip .3cm
\noindent {\bf Definition of $\krv$.} Let $g,h$ be Lie polynomials of 
homogeneous degree $n\ge 3$, and let $D_{g,h}$ be the tangential derivation
of ${\rm Lie}[x,y]$ defined by $y\mapsto [y,g]$, $x\mapsto [x,h]$. The
derivation $D_{g,h}$ lies in $\krv$ if it is a {\it special derivation},
i.e.~if it also satisfies $D_{g,h}(x+y)=0$,
and writing $g=g_xx+g_yy$ and $h=h_xx+h_yy$, it further satisfies the 
{\it divergence condition} for some constant $c$:
$$tr(h_xx+g_yy)=nc\,tr((x+y)^n-x^n-y^n)\eqno(1)$$
in the space $Tr\langle x,y\rangle$ which is the quotient of the
polynomial ring ${\Bbb Q}\langle x,y\rangle$ by cyclic permutations
of monomials. The space of derivations $\krv$ spanned by the
$D_{g,h}$ is a Lie algebra under the usual
bracket of derivations (see [1], just before Prop. 4.1).
\vskip .4cm
\noindent Calculation in low weights and comparison with the known injective
Lie algebra morphisms from the Grothendieck-Teichm\"uller Lie algebra $\grt$ 
into $\ds$ (Furusho) and into $\krv$ (Alekseev-Torossian) 
has long suggested that for each $f\in \ds$, setting $g=f(-x-y,-y)$,
there exists an element $h\in {\rm Lie}[x,y]$ such that
$f\mapsto D_{g,h}$
provides an injective Lie algebra morphism from $\ds$ to $\krv$. 
A proof of this fact was given in
[12], but the argument was partly based on a stated result by J.~\'Ecalle, 
which turned out to have an insufficient proof in the literature, and which no one has been able to prove since, so that the result is considered 
still an open question.
\vskip .3cm
\noindent {\bf Theorem 1.} {\it Let $f\in \ds$ and let $g=f(-x-y,-y)$. 
Then there exists a Lie polynomial $h\in {\rm Lie}[x,y]$ such that
the derivation $D_{g,h}$ of ${\rm Lie}[x,y]$ defined by $x\mapsto [x,h]$
and $y\mapsto [y,g]$ also satisfies 
\vskip .2cm
\noindent (i) $D_{g,h}(x+y)=0$, i.e.~$D_{g,h}$ is a special derivation, and
\vskip .1cm
\noindent (ii) $D_{g,h}$ satisfies the divergence condition (1).
\vskip .2cm
\noindent The map $f\mapsto D_{g,h}$ gives an injective Lie algebra 
morphism
$$\ds\hookrightarrow \krv.$$}\par
\noindent {\bf Remark.} This morphism is compatible with the injective Lie
morphism from the Grothendieck-Teichm\"uller Lie algebra $\grt$ into $\ds$, which 
maps $f(x,y)\mapsto f(x,-y)$, proven by Furusho in [7], and the injective Lie
morphism from $\grt\hookrightarrow \krv$ proven by Alekseev and Torossian in [1],
which maps $f(x,y)$ to the derivation of ${\rm Lie}[x,y]$ taking $x\mapsto [x,f(-x-y,x)]$
and $y\mapsto [y,f(-x-y,y)]$. Thus Theorem 1 gives a commutative diagram of 
injective maps
$$\xymatrix{\grt\ar[dr]\ar[rr]&&\krv\\
&\ds\ar[ur].}$$
\vskip .2cm
\noindent The proof we give here is based on Theorem 2 below, along with
technical Lemmas 3-6 and Corollaries 1 and 2. Theorem 2 is a compendium 
of known results from mould theory, originally proven in various
articles, but all of which can be found in in the basic mould
theory reference work [11].
\vskip .3cm
\noindent Following the notation for ${\rm Lie}[x,y]$, let ${\rm Lie}[a,b]$ 
denote the {\it degree-completed} free Lie algebra on two variables $a$
and $b$. Let
$$t_{01}=Ber_b(-a),\ \ \  t_{02}=Ber_{-b}(a),\ \ \ t_{12}=[a,b],\eqno(2)$$
be three elements of ${\rm Lie}[a,b]$, 
where 
$$Ber_x(y)={{ad(x)}\over{e^{ad(x)}-1}}(y)=\sum_{n\ge 0} {{B_n}\over{n!}}ad(x)^n(y),$$
where $B_n$ denotes the $n$-th Bernoulli number.
They satisfy
$t_{01}+t_{02}+t_{12}=0$. These elements generate a free Lie subalgebra
${\rm Lie}[t_{01},t_{02}]$ on two generators inside ${\rm Lie}[a,b]$, 
so we have a injective Lie algebra morphism
$$\eqalign{{\rm Lie}[x,y]&\rightarrow {\rm Lie}[a,b]\cr
x&\mapsto t_{01}\cr
y&\mapsto t_{02}}\eqno(3)$$
which also maps $z:=-x-y$ to $t_{12}$.
\vskip .3cm
\noindent The {\it push-operator} is defined on monomials in $x$ and $y$
(resp,~in $a$ and $b$) as a cyclic permutation of the power of $x$'s (resp.~$a$'s) between the $y$'s (resp.~$b$'s):
$$\cases{push(x^{i_0}yx^{i_1}y\cdots yx^{i_r})=x^{i_r}yx^{i_0}y\cdots yx^{i_{r-1}}\cr
push(a^{i_0}ba^{i_1}b\cdots ba^{i_r})=a^{i_r}ba^{i_0}b\cdots ba^{i_{r-1}},}$$
and extended to polynomials by linearity. A polynomial is said to be 
{\it push-invariant} if it is fixed by the push-operator.
\vskip .3cm
\noindent {\bf Theorem 2.} {\it Let $f(x,y)\in \ds$, and assume that
$f$ is homogeneous of degree $n\ge 3$, with (minimal) depth $d\ge 1$. 
Then there exists a unique derivation $D$ of ${\rm Lie}[a,b]$ satisfying
the following properties:
\vskip .2cm\noindent
(1) $D(a)$ is a push-invariant element of ${\rm Lie}[a,b]$, 
\vskip .1cm\noindent
(2) $D([a,b])=0$,
\vskip .1cm\noindent
(3) $D(t_{02})=[f(t_{02},-t_{12}),t_{02}]$, 
\vskip .1cm\noindent
(4) $D(a)$ has only terms of odd degree in $a,b$.}
\vskip .2cm
\noindent Proof. These results, originally proven in [13],
can all be found in \S 5 of the mould theory reference text [11].
More precisely, the existence of a derivation $D$ of
${\rm Lie}[a,b]$ satisfying (1) and (2) is proven in Theorem 5.6.1 
the action on $t_{02}$ is proved in Theorem 5.6.1 (ii),
the uniqueness of $D$ in Theorem 5.6.1 (iii), and property
(4) in Theorem 5.6.1 (iv).\hfill{$\square$}
\vskip .3cm
\noindent {\bf Corollary 1.} {\it Let $D$ be as in Theorem 2 (for 
$f\in \ds$).  Then we have
$$\cases{D(t_{01})=[f(t_{01},-t_{12}),t_{01}]\cr
D(t_{02})=[f(t_{02},-t_{12}),t_{02}]\cr
D(t_{12})=0.}\eqno(4)$$}
\vskip .1cm
\noindent Proof.  Let $\iota$ denote the 
involutive automorphism of ${\rm Lie}[a,b]$ defined by 
$$\iota(a)=-a,\ \ \ \iota(b)=-b.$$
We claim that $D$ commutes with $\iota$ on ${\rm Lie}[a,b]$.
To check this, we consider the derivation $D'=\iota\circ D\circ \iota$
of ${\rm Lie}[a,b]$, and compare $D'$ with $D$ on $a$ and $[a,b]$. 
On $a$, we find that
$$D'(a)=(\iota\circ D\circ \iota)(a)=
\iota\circ D(-a)=-\iota\bigl(D(a)\bigr).$$
But since $D(a)$ has only odd-degree terms by Theorem 2 (iv), we have
$\iota\bigl(D(a)\bigr)=-D(a)$ and therefore $D'(a)=D(a)$.
On $[a,b]$, since $\iota([a,b])=[a,b]$, we have
$$D'([a,b])=\iota\circ D([a,b])=0=D([a,b]).$$
Therefore $D$ and $D'$ agree on $a$ and $[a,b]$, and since a derivation
annihilating $[a,b]$ is uniquely determined by its value on $a$,
we have $D'=D$, proving that $D$ commutes with $\iota$. 

The second and third lines of (4) are proved in Theorem 2 (ii) and (iii).
To prove the first line of (4), we simply observe that 
$t_{01}=\iota(t_{02})$, so
$$D(t_{01})=D\bigl(\iota(t_{02})\bigr)=\iota\bigl(D(t_{02})\bigr)
=\iota\bigl([f(t_{02},-t_{12}),t_{02}]\bigr)
=[f(t_{01},-t_{12}),t_{01}].$$
This concludes the proof.\hfill{$\square$}
\vskip .3cm
\noindent For any $g\in {\rm Lie}[x,y]$, write $g=g_xx+g_yy$ and define the {\it partner}
$g'$ of $g$ by the formula
$$g'=\sum_{i\ge 0} {{(-1)^{i-1}}\over{i!}}x^iy\partial^i_x(g_x)\eqno(5)$$
where $\partial_x$ is the derivation defined by $\partial_x(x)=1$,
$\partial_x(y)=0$.  It is shown in Theorem 2.1 of [12] that if $D_{g,h}$ denotes
the tangential derivation of ${\rm Lie}[x,y]$ defined by 
$y\mapsto [y,g]$, $x\mapsto [x,h]$ with $g,h\in {\rm Lie}[x,y]$, then 
$D_{g,h}(x+y)=0$ if and only if $g$ is push-invariant and $h=-g'$. 
Observe that if $E_{g,g'}$ denotes the derivation of ${\rm Lie}[x,y]$ given by 
$x\mapsto g$, $y\mapsto g'$, and if $g$ is push-invariant, then
$$E_{g,g'}([x,y])=[g,y]+[x,g']=0.\eqno(6)$$
\vskip .1cm
\noindent {\bf Corollary 2.} {\it If $f\in \ds$ then the element
$f(y,-z)$ is push-invariant as a polynomial in $x$ and $y$.}
\vskip .2cm
\noindent Proof. 
Since we have the Lie algebra isomorphism (3), the derivation of
${\rm Lie}[t_{01},t_{02}]$ given in (4) pulls back to a derivation of
${\rm Lie}[x,y]$ given by
$$\cases{x\mapsto [f(x,-z),x]=[x,-f(x,-z)]\cr
y\mapsto [f(y,-z),y]=[y,-f(y,-z)]\cr
x+y\mapsto 0.}\eqno(7)$$ 
Since $x+y\mapsto 0$, we thus see that $-f(y,-z)$ must be 
push-invariant as a polynomial in $x$ and $y$,
so $f(y,-z)$ is as well.\hfill{$\square$} 
\vskip .3cm \noindent 
\noindent {\bf Remark.} Since (7) shows that an element $f\in \ds$
induces a special derivation of ${\rm Lie}[x,y]$, one might think that
this could provide a morphism from $\ds$ into $\krv$. However,
we were not able to prove that the derivation in (7) satisfies the 
divergence condition (although this seems very probable). We will 
now show that a slightly modified version of this derivation works.
\vskip .2cm
\noindent {\bf Lemma 4.} {\it Let $g\in {\rm Lie}[x,y]$ be a push-invariant
polynomial, and let $\epsilon\in {\rm Der}{\rm Lie}[x,y]$ be the derivation defined
by $\epsilon(x)=0$, $\epsilon(y)=x$. Then $\epsilon(g)$ is push-invariant.}
\vskip .2cm
\noindent Proof. The fact that $g$ is push-invariant means that for every
monomial $x^{i_0}y\cdots yx^{i_r}$ appearing in $g$, the push-images of this
monomial appear with the same coefficient; thus all the terms
in the orbit under the push-powers of the push-operator appear in $g$
with the same coefficients; in other words $g$ is a sum of push-orbits
$$x^{i_0}y\cdots yx^{i_r}+x^{i_r}yx^{i_0}y\cdots yx^{i_{r-1}}+\cdots+
x^{i_1}y\cdots yx^{i_r}yx^{i_0}\eqno(8)$$
with various coefficients. The image of a push-orbit polynomial under
$\epsilon$ consists in the sum of the $r(r+1)$ (not necessarily distinct) 
monomials obtained by changing exactly one $y$ to an $x$ in each of the $r+1$ 
monomials in (8). This set of $r(r+1)$ monomials breaks up into $r+1$
push-orbits, one for each pair of indices $(i_j,i_{j+1})$ and $(i_r,i_0)$,
by changing the $y$ between $i_j$ and $i_{j+1}$ to an $x$ in each of the
$r$ monomials of (8) where such a $y$ exists (i.e.~every monomial of
(8) except for $x^{i_{j+1}}y\cdots yx^{i_j}$). 
Thus the image of (8) by $\epsilon$ is the sum of $r+1$ push-orbits of
length $r$, so $\epsilon(g)$ is push-invariant.\hfill{$\square$}
\vskip .3cm
\noindent {\bf Lemma 5.} {\it Let $f\in \ds$. Then $f(z,-y)$ is push-invariant
as a polynomial in $x$ and $y$.}
\vskip .2cm
\noindent Proof.  Set $g(x,y)=f(y,-z)$. By Corollary 2, we know this element
is push-invariant. Applying the automorphism $x\mapsto x$, $y\mapsto z$,
$z\mapsto y$ (with $z=-x-y$) translates this equality into $g(x,z)=f(z,-y)$.
So the Lemma comes down to proving that if $g(x,y)$ is push-invariant then
$g(x,z)$ is also push-invariant.  For $1\le r\le d$, let $g^r$ denote the 
depth $r$ part of $g$ (i.e. the part of homogeneous $y$-degree $r$). 
Then the depth $r$ part of $g(x,z)$ is given by
$$(-1)^r\Bigl(g^r+\epsilon(g^{r+1})+{{1}\over{2}}\epsilon^2(g^{r+2})+\cdots+{{1}\over{(d-r)!}}\epsilon^{d-r}(g^d) \Bigr).$$
But since push-invariance is a depth-by-depth condition, $g^r$ is push-invariant
for all $1\le r\le d$, and thus by Lemma 4 so are the terms 
$\epsilon^i(g^r)$ for all $i$ and $r$.  Thus $g(x,z)=f(z,-y)$ is push-invariant.\hfill{$\square$}
\vskip .3cm
\noindent Before proving Theorem 1, we need one more definition.
\vskip .3cm
\noindent {\bf Definition.}
For any monomial $u$ in $x$ and $y$, let $r$ denote the number
of $y$'s in $u$, and let $Push(u)$ denote the list of $r+1$ (possibly 
not distinct) monomials obtained from $u$ by applying the push-operator $r+1$ times. 
For $n\ge 3$, a polynomial $p\in {\Bbb Q}\langle x,y\rangle$ of homogeneous degree 
$n-1$ is said to be {\it push-constant} if $(p|y^{n-1})=0$ and for all monomials $u$
of degree $n-1$ there exists a constant $c$, which is equal to $0$ if $n-1$ is odd, 
such that
$$\sum_{u'\in Push(u)} (p\,|\,u')=c.\eqno(9)$$
\vskip .1cm
\noindent {\bf Proof of Theorem 1.} 
Let $f\in \ds$ be an element of homogeneous degree
$n\ge 3$. By Lemma 5, 
$g=f(z,-y)$ is push-invariant. Thus by Theorem 2.1 of [12], $g$ admits a 
unique partner $h\in {\rm Lie}[x,y]$ (defined by (5)) such that $[x,h]+[y,g]=0$.
Let $D_{g,h}$ denote the special derivation defined by
$$\cases{D_{g,h}(y)=[y,g]=[y,f(z,-y)]\cr
D_{g,h}(x)=[x,h]\cr
D_{g,h}(x+y)=0.}\eqno(10)$$
We will show that $D_{g,h}$ lies in $\krv$, 
i.e.~that it satisfies the divergence condition. We will
derive this from the following known result on $f\in \ds$ (cf.~Corollary 1 of [3],
also stated in Theorem 3.4 of [12]). Let $c=(f|x^{n-1}y)$ and write
$f=f_xx+f_yy$. Then $(f_y|y^{n-1})=0$ and for every degree $n-1$ monomial 
$u\ne y^{n-1}$, letting $r$ denote the number of $y$'s in $u$, 
$f_y$ satisfies the identity
$$\sum_{u'\in Push(u)} (f_y\,|\,u')=(-1)^{r-1}c.\eqno(11)$$
A well-known result states that an element $f\in \ds$ of homogeneous degree $n$ has 
(minimal) $y$-degree of the same parity as $n$ (cf.~for example Corollary 3.4.4 of [11]),
so if $n$ is even, the constant $c=0$.
Setting $\tilde f(x,y)=f(x,-y)$, we immediately deduce from (11) that 
$\tilde f_y$ is push-constant.
\vskip .2cm
\noindent It is shown in Lemma 3.5 of [12] that if $\tilde f\in {\rm Lie}[x,y]$ is
such that $\tilde f_y$ is push-constant for $c$ (with $c=0$ if $n$ is even), 
and if we set $g=\tilde f(z,y)$, then $g_y-g_x$ is push-constant for the same 
value of $c$. A fortiori this holds for our choice of $f\in \ds$.
\vskip .2cm
\noindent Let us use the notation $\overleftarrow{u}$ for a monomial $u$ written 
backwards, and the same notation for a polynomial. Any Lie polynomial of
homogeneous degree $n$ satisfies $g=(-1)^{n-1}\overleftarrow{g}$. Writing 
$g=yg^y+xg^x=g_xx+g_yy$, we see that
$$g=(-1)^{n-1}\overleftarrow{g}=(-1)^{n-1}(\overleftarrow{g_xx+g_yy})=(-1)^{n-1}x\overleftarrow{g_x}+(-1)^{n-1}y\overleftarrow{g_x}=xg^x+yg^y,$$
so $g^y=(-1)^{n-1}\overleftarrow{g_y}$, $g^x=(-1)^{n-1}\overleftarrow{g_x}$ and
$$g^y-g^x=(-1)^{n-1}\Bigl(\overleftarrow{g_y-g_x}\Bigr).$$
If a polynomial $p$ is push-constant then the backwards-written polynomial
$\overleftarrow{p}$ is also push-constant (and multiplying $p$ by any
constant also preserves push-constance), so we find that $g^y-g^x$ is
push-constant. It is shown in [12] (between (3.3) and (3.5)) that 
$g^y-g^x$ being push-constant for a constant $c$ is equivalent to 
$$tr\bigl((g^y-g^x)y\bigr)=c\,tr\bigl((x+y)^n-x^n-y^n\bigr).\eqno(12)$$
Now, writing $g=xg^x+yg^y=g_xx+g_yy$ and $h=xh^x+yh^y=h_xx+h_yy$, we can
rewrite the equality $[x,h]+[y,g]=0$ as $[x,h]=[g,y]$ and expand it as
$$xh_xx+xh_yy-xh^xx-yh^yx=xg^xy+yg^yy-yg_xx-yg_yy,$$
which by comparing words that start and end with the same letter shows that
$$h_x=h^x,\ \ h_y=g^x,\ \ h^y=g_x,\ \  g^y=g_y.$$
Thus the left-hand side of (12) can be rewritten
$tr\bigl((g_y-h_y)y\bigr)$, but since the trace of a Lie element is
always zero, we have $tr(-h_yy)=tr(h_xx)$, so (12) can be rewritten as
$$tr\bigl(g_yy+h_xx\bigr)=c\,tr\bigl((x+y)^n-x^n-y^n\bigr),\eqno(13)$$
which is exactly the divergence condition (1). This
shows that the special derivation $D_{g,h}$ associated to $f\in \ds$ by (7) satisfies 
the divergence relation, and therefore $f\mapsto D_{g,h}$ gives an injective 
linear map from $\ds$ to $\krv$. 
\vskip .2cm
\noindent We complete the proof by showing that this map is a Lie algebra morphism. Recall that 
the Lie bracket on $\ds$ is the Poisson bracket given by
$$\{f,f'\}-D_f(f')-D_{f'}(f)+[f,f'].\eqno(14)$$
Consider the automorphism $\alpha$ of ${\rm Lie}[x,y]$ defined by $\alpha(x)=x$,
$\alpha(y)=-y$. Let $\tilde f=\alpha(f)$, and let $\widetilde\ds=\alpha(\ds)$.
We have
$$D_{\alpha(f)}(\alpha(f'))-D_{\alpha(f')}(\alpha(f))+[\alpha(f),\alpha(f')]=
\alpha\bigl(D_f(f')-D_{f'}(f)+[f,f']\bigr)=\alpha(\{f,f'\}),$$
so the Poisson bracket in (14) also defines a Lie bracket on $\ds'$, and
$\alpha:\ds\rightarrow\ds'$ is a Lie algebra isomorphism. 
is an injective Lie algebra morphism.
\vskip .2cm
\noindent Let $\beta$ denote the involutive automorphism of ${\rm Der}{\rm Lie}[x,y]$ defined by
$\beta(z)=x$ and $\beta(y)=y$. Conjugation with $\beta$ defines an automorphism of
the space ${\rm Der}{\rm Lie}[x,y]$. Thus the following composition of maps
forms an injective Lie algebra from $\ds$ to ${\rm Der}{\rm Lie}[x,y]$:
$$\ds\rightarrow\ds'\rightarrow {\rm Der}{\rm Lie}[x,y]\rightarrow {\rm Der}{\rm Lie}[x,y]$$
$$\ \ \ \ \ f\mapsto \alpha(f)\mapsto \ D_{\alpha(f)}\ \ \mapsto \ \ \ \beta^{-1}\circ D_{\alpha(f)}\circ \beta.$$
But this composition of maps takes $f\in \ds$ to the derivation mapping
$z\mapsto 0$ and $y\mapsto [y,f(z,-y)]$, which is precisely $D_{g,h}$.
This shows that the map $f\mapsto D_{g,h}$ is an injective Lie algebra morphism 
from $\ds$ to $\krv\subset {\rm Der}{\rm Lie}[x,y]$, 
completing the proof of Theorem 1.\hfill{$\square$}

\vskip 1cm
\noindent {\bf References.} 
\vskip .4cm
\noindent [1] A.~Alekseev, C.~Torossian, The Kashiwara-Vergne conjecture and Drinfeld's associators, {\it Annals of Math.} {\bf 175} (2012), 415-463.
\vskip .2cm
\noindent [2] S.~Baumard, Thesis, 2014.
\vskip .2cm
\noindent [3] S.~Carr, L.~Schneps, Combinatorics of the double shuffle Lie algebra, in Galois-Teichmüller theory and Arithmetic Geometry, H. Nakamura, F. Pop, L. Schneps, A. Tamagawa, eds., Adv. Stud. Pure Math. 63, Mathematical Society of Japan, 2012. 
\vskip .2cm
\noindent [4] D.~Dorigoni et al, Canonicalizing zeta generators: genus zero and genus one, arXiv: 2406.05099v2.
\vskip .2cm
\noindent [5] J.~\'Ecalle, The flexion structure of dimorphy: flexion units, singulators, generators, and the enumeration of multizeta irreducibles, in {\it Asymptotics in Dynamics, Geometry and PDEs; Generalized Borel Summation II}, O. Costin, F. Fauvet, 
F. Menous, D. Sauzin, eds., Edizioni della Normale, Pisa, 2011.
\vskip .2cm
\noindent [6] B.~Enriquez, H.~Furusho, Double shuffle Lie algebra
and special derivations, preprint 2025.
\vskip .2cm
\noindent [7] H.~Furusho, Double shuffle relation for associators, {\it Annals of Math.}
 {\bf 174} (2011), 341-360.
\vskip .2cm
\noindent [8] G.~Racinet, Doubles m\'elanges des polylogarithmes 
multiples aux racines de l'unit\'e. {\it Publ. Math. Inst.
Hautes Etudes Sci.}, {\bf 95} (2002), 185–231.
\vskip .2cm
\noindent [9] E.~Raphael, L.~Schneps, On linearised and elliptic versions of the Kashiwara-Vergne Lie algebra, arXiv:1706.08299. 
\vskip .2cm
\noindent [10] A.~Salerno, L.~Schneps, Mould theory and the double shuffle Lie algebra 
structure, in {\it Periods in Quantum Field Theory and Arithmetic}, 
J.~Burgos Gil, K.~Ebrahimi-Fard, H.~Gangl, eds.,
Springer PROMS 314, 2020, 399-430.
\vskip .2cm
\noindent [11] L.~Schneps, ARI, GARI, Zig and Zag, arXiv:1507.01534
\vskip .2cm
\noindent [12] L.~Schneps, Double shuffle and Kashiwara-Vergne Lie algebras,
{\it J. Alg} {\bf 367} (2012), 54-74.
\vskip .2cm
\noindent [13] L.~Schneps, Elliptic multizeta values, Grothendieck-Teichm\"uller and mould theory, Ann. Math. Qu\'ebec {\bf 44} (2) (2020), 261-289,
updated version with minor corrections arXiv:1506.09050.
\vskip .2cm

\bye